\documentclass [11pt,twoside,a4paper]{article}
\usepackage{amsfonts}
\usepackage{amsmath}
\usepackage{amstext}
\usepackage{amssymb}
\usepackage{mathrsfs}
\usepackage{amscd}
\usepackage{xypic}
\usepackage{epsf}              
\usepackage{graphicx}          
\usepackage{fancybox}          
\usepackage{color}             
\usepackage{fancyhdr}
\usepackage[hang,footnotesize]{caption2}  
\usepackage{mathbbold}          

\def\mathcal{\mathscr}
\newfont{\aaa}{cmb10 at 19pt}
\newfont{\bbb}{cmb10 at 11pt}
\newtheorem{thm}{Theorem}[section]
\newtheorem{cor}[thm]{Corollary}
\newtheorem{lem}[thm]{Lemma}
\newtheorem{prop}[thm]{Proposition}

\newtheorem{rem}[thm]{Remark}
\newtheorem{exm}[thm]{Example}
\makeatletter
\def\@evenhead{
\vbox{\hbox to \textwidth {}{\hspace{0mm}{\footnotesize
\thepage}}{\hspace{10cm} {\footnotesize Mu-Fa CHEN et al.}} \protect\vspace{1truemm}\relax \hrule depth0pt
height0.15truemm width\textwidth}}
\def\@evenfoot{}
\def\@oddhead{\vbox{\hbox to \textwidth
{{\hspace{0cm}{\footnotesize Mixed eigenvalues of $p\,$-Laplacian}\hfill{\footnotesize
\thepage}}\hspace{0mm}}{} \protect\vspace{1truemm}\relax\hrule
depth0pt height0.15truemm width\textwidth}}
\def\@oddfoot{}
\makeatother

\def\d{{\rm d}}

\def\scr{\mathscr}
\def\supp{\text{\scriptsize\rm supp\,}}

\newcommand{\rf}[2]{[\ref{#1}; #2]}

\begin{document}

\thispagestyle{empty} \thispagestyle{fancy} {
\fancyhead[lO,RE]{\footnotesize  Front. Math. China \\
DOI 10.1007/s11464-015-0375-0\\[3mm]
}
\fancyhead[RO,LE]{\scriptsize \bf 
} \fancyfoot[CE,CO]{}}
\renewcommand{\headrulewidth}{0pt}


\setcounter{page}{1}
\qquad\\[8mm]

\noindent{\aaa{Mixed  Eigenvalues of  {\LARGE $\pmb p$\,}-Laplacian}}\\[1mm]

\noindent{\bbb Mu-Fa CHEN$^1$,\quad Ling-Di WANG$^{1,2}$,\quad Yu-Hui ZHANG$^1$}\\[-1mm]

\noindent\footnotesize{1 Beijing Normal University, Beijing 100875, China\\
2 Henan University, Kaifeng, Henan 475004, China}\\[6mm]

\vskip-2mm \noindent{\footnotesize$\copyright$ Higher Education
Press and Springer-Verlag Berlin Heidelberg 2013} \vskip 4mm

\normalsize\noindent{\bbb Abstract}\quad The mixed principal eigenvalue of
$p\,$-Laplacian (equivalently, the optimal constant of weighted
Hardy inequality in $L^p$ space) is studied in this paper. Several
variational formulas for the eigenvalue are presented. As
applications of the formulas, a criterion for the positivity of the
eigenvalue is obtained. Furthermore,  an
approximating procedure and some explicit estimates are presented case by case. An example is
included to illustrate the power of the results of the paper.\vspace{0.3cm}

\footnotetext{Received August 13, 2013; accepted March 3,
2014\\
\hspace*{5.8mm}Corresponding author: Ling-Di WANG, E-mail:wanglingdi@henu.edu.cn}

\noindent{\bbb Keywords}\quad $p\,$-Laplacian, Hardy inequality in $L^p$ space, mixed boundaries,
explicit estimates, eigenvalue, approximating procedure\\
{\bbb MSC}\quad 60J60, 34L15\\[0.4cm]

\setcounter{section}{1}
\setcounter{thm}{0}

\noindent{\bbb 1\quad Introduction}\\[0.1cm]
\noindent
 As a natural extension of Laplacian from linear to nonlinear, $p\,$-Laplacian plays a typical role
     in mathematics,
     especially in nonlinear analysis. Refer to \cite{CDF2004,LE2011}  for recent progresses on this subject.
     Motivated by the study on stability speed, we come to this topic, see \cite{Chen2005,Chen2010} and
     references therein.
     The present paper is a continuation of \cite{CWZ2013-1} in which the  estimates of the mixed
     principal eigenvalue
     for discrete $p\,$-Laplacian were carefully studied. This paper deals with the same problem but for
     continuous $p\,$-Laplacian,
     its principal eigenvalue is equivalent to the optimal constant in the weighted Hardy inequality.
     Even though the discrete case is often harder than the
     continuous one, the latter has its own difficulty. For
     instance, the existence of the eigenfunction is rather hard in
     the nonlinear context, but it is  not a problem in the
     discrete situation.
     Similar to the case of $p=2$ (\cite{Chen2010,CWZ2013}), there are four types of boundaries:
     Neumann (denoted by code ``N'') or Dirichlet (denoted by code ``D'') boundary at the left- or
     right-endpoint
     of the half line $[0, D]$.  In \cite{JM2012}, Jin and Mao studied a class  of weighted Hardy
     inequality and
     presented two variational formulas in the DN-case. Here, we study  ND-case carefully and
     add some results
     to \cite{JM2012}. The DD- and NN-cases will be handled elsewhere.
     Comparing with our previous study, here the general weights are allowed.

The paper is organized as follows. In the next section, restricted
in the ND-case, we introduce the main results: variational formulas
and the  basic estimates for the optimal constant (cf. \cite{KMP2006,KMP2007}). As an application, we improve the basic
estimates step by step through an approximating procedure. To
illustrate the power of the results, an example is included. The
sketched proofs of the results in Section 2 are presented in Section
3. For another mixed case: DN-case studied in \cite{JM2012}, some
complementary are presented in Section 4.
\noindent\\[4mm]

\setcounter{section}{2}
\setcounter{thm}{0}

\noindent{\bbb 2\quad ND-case}\\[0.1cm]
\noindent
 Let $\mu$, $\nu$ be two positive Borel measures on $[0, D]$,
$D\leqslant \infty$ (replace $[0,D]$ by $[0, D)$ if $D=\infty$), $\d\mu=u(x)\d x$ and $\d\nu=v(x)\d x$.
Next, let
$$L_pf=\big(v |f'|^{p-2}f'\big)',\qquad p>1.$$
Then the eigenvalue problem with ND-boundary conditions reads:
\begin{equation}\begin{cases}
\text{Eigenequation}: &L_pg(x)=-\lambda u(x)|g|^{p-2}g(x);\\
\text{ND-boundaries}: &g'(0)=0,\quad g(D)=0\text{ if } D<\infty.
\label{CND-eigen}\end{cases}
\end{equation}
If $(\lambda, g)$ is a solution to the eigenvalue problem above, $g\neq0$,
then we call $\lambda$  an  `eigenvalue' and $g$ is an
`eigenfunction' of $\lambda$. When $p=2$, the operator $L_p$
defined above returns  to the diffusion operator defined in
\cite{CWZ2013}: $u^{-1}(vf')'$,
 where $u(x)\d x$ is the  invariant measure of the diffusion process and $v$ is a Borel
 measurable function related to its recurrence criterion.
 For $\alpha\leqslant\beta$, define
$$\mathscr{C}[\alpha, \beta]=\big\{f: f \text{ is  continuous on }[\alpha, \beta] \big\},$$
  $$\mathscr{C}^k(\alpha, \beta)=\big\{f: f \text{ has continuous derivatives of order } k
   \text{ on } (\alpha, \beta)\big\},\quad k\geqslant1,$$
  and
$$\mu_{\alpha, \beta}(f)=\int_{\alpha}^{\beta}f\d \mu,\quad\quad D_p^{\alpha, \beta}(f)
=\int_{\alpha}^{\beta}|f'|^p\d\nu.$$
Similarly, one may define $\scr{C}(\alpha, \beta)$. In this section, we  study the first eigenvalue (the minimal one), denoted by $\lambda_p$,
described by the following classical variational formula:
\begin{equation}\label{N-D-eif1}
\lambda_p=\inf\big\{{D_p(f)}: f\in \mathscr{C}_K[0, D], \mu(|f|^p)=1, f(D)=0\text{ if } D<\infty\big\},
\end{equation}
where $\mu(f)=\mu_{0,D}(f)$, $D_p(f)=D_p^{0,D}(f)$ and
$$\mathscr{C}_K[\alpha, \beta]=\big\{f\in\scr{C}[\alpha, \beta]: v^{p^*-1}f'\in\scr{C}(\alpha, \beta)
\text{ and } f \text{ has compact support}\big\},$$
with $p^*$ the conjugate number of $p$ (i.e., $p^{-1}+{p^*}^{-1}=1$).
When $p=2$, it reduces to the linear case studied in \cite{CWZ2013}.
Thus, the aim of the paper is extending the results in linear case
($p=2$) to nonlinear one. Set
 $$\mathscr{A}[\alpha, \beta]=\big\{f: f \text{ is absolutely continuous on }[\alpha, \beta]\big\}.$$
 As will be proved soon  (see Lemmas \ref{CNL0} and \ref{CNL01}), we can rewrite  $\lambda_p$ as
\begin{equation}\label{CNDf1}\tilde\lambda_{*,p}:=\inf\big\{D_p(f): \mu(|f|^p)=1,  f\in\mathscr{A}[0, D],
 f(D)=0\big\}.\end{equation}

By making inner product   with $g$ on  both sides of  eigenequation \eqref{CND-eigen} with respect
to the Lebesgue measure over $(\alpha, \beta)$, we obtain
$$\lambda\mu_{\alpha, \beta}(|g|^p)=D_p^{\alpha, \beta}(g)-\big(v|g'|^{p-2}g'g\big)\big|_{\alpha}^{\beta}.$$
 Moreover, since $g'(0)=0$, we have
$$\lambda\mu(|g|^p)=D_p(g)-\big(v|g'|^{p-2}g'g\big)(D),$$
 where, throughout this paper, $f(D):=\lim_{x\to D} f(x)$ provided $D=\infty$. Hence, with
$$\mathscr{D}(D_p)=\{f: f \in \mathscr{A}[0, D],\; D_p(f)<\infty\},$$
$A:=\lambda_p^{-1}$ is the optimal constant  of the following {\it weighted Hardy inequality}:
$$\aligned
&\text{Hardy inequality}:\quad \mu(|f|^p)\leqslant A D_p(f),\quad f\in \mathscr{D}(D_p);\\
&\text{Boundary condition}:\quad f(D)=0.
\endaligned$$
Note that the boundary condition ``$f'(0)=0$'' is unnecessary in the
inequality.

Throughout this paper, we concentrate on $p\in (1, \infty)$ since
the degenerated cases that either $p=1$ or $\infty$ are often easier
to handle (cf. \rf{BA1990}{Lemmas 5.4 and 5.6 on pages 49 and 56,
respectively}).

\vspace{2mm}
\noindent{\bbb Main notation and results}

\vspace{1mm}\noindent
For $p>1$, let  $p^*$ be its conjugate number.
Define $\hat{v}(x)=v^{1-p^*}(x)$ and $\hat\nu(\d x)=\hat{v}(x)\d x$.
We use the following hypothesis throughout the paper:
$$\aligned u, {\hat v}\; \text{\it are locally integrable with respect to the Lebesgue measure on } [0, D],
\endaligned$$
without mentioned time by time.

Our main operators are defined as follows.
$$\aligned
I(f)(x)&=-\frac{1}{\big(v f'|f'|^{p-2}\big)(x)}\int_{0}^{x}f^{p-1}\d\mu
\qquad(\text{single integral form}),
\\ I\!I(f)(x)&=\frac{1}{f^{p-1}(x)}\bigg[\int_{(x, D)\cap{\supp}(f)}\!\!\!\hat{v}(s)
\bigg(\int_{0}^{s}f^{p-1}\d\mu\bigg)^{p^*-1}\d s\bigg]^{p-1}\\
&\hskip7cm\quad(\text{double integral form}),\\
R(h)(x)&\!=\!u(x)^{-1}\big[-|h|^{p-2}\big(v'h\!+\!(p\!-\!1)(h^2\!\!+h')v\big)\big](x)\;
\quad(\text{differential
form}).
\endaligned$$
These operators have domains, respectively, as follows.
$$\aligned&\mathscr{F}_{I}=\{f\in\mathscr{C}[0,D]:v^{p^*-1}f'\in \mathscr{C}(0,D),  f|_{(0,D)}>0,  f'|_{(0, D)}<0\},\\
&\mathscr{F}_{I\!I}=\{f:f\in\mathscr{C}[0,D], f|_{(0, D)}>0\},\\
&\mathscr{H}\!=\!\{h: h\in\mathscr{C}^{1}(0,
D)\cap\mathscr{C}[0,D], h(0)=0, h|_{(0,D)}<0 \text{ if }\hat{\nu}(0,D)<\infty,\\
 &\hskip 1.6cm\text{ and } h|_{(0,D)}\leqslant 0 \text{ if }\hat{\nu}(0,D)=\infty\},
\endaligned$$
where $\nu(\alpha,\beta)=\int_{\alpha}^{\beta}\d\nu$ for a measure
$\nu$. To avoid the non-integrability problem, some modifications of
these sets are needed for studying the upper estimates.
$$\aligned
&\widetilde{\mathscr{F}}_{I}=\big\{f\in\mathscr{C}[x_0, x_1]: v^{p^*-1}f'\in\mathscr{C}(x_0, x_1), f'|_{(x_0, x_1)}<0 \text{ for some }\\
&\hskip 1.8cm x_{0}, x_{1}\in(0, D)\text{ with } x_0<x_1, \text{ and
}f=f(\cdot\vee x_{0}) \mathbbold{1}_{[0, x_{1})}\big\}.
\endaligned$$
$$\aligned
\widetilde{\mathscr{F}}_{I\!I}=\big\{f:f=f \mathbbold{1}_{[0, x_{0})}\text{ for some }\ x_{0}\in(0,
D) \text{ and } f\in\mathscr{ C}[0, x_{0}] \big\},\endaligned$$
$$\aligned\widetilde{\mathscr{H}}=\Big\{h: \;&\exists x_{0}\in(0,
D)\text{ such that } h\in\mathscr{C}[0,x_{0}]\cap\mathscr{C}^{1}(0,
x_{0}),\; h|_{(0, x_{0})}<0,\\& h|_{[x_0, D]}=0,  h(0)=0, \text{ and
}\sup_{(0, x_0)}\big(v'h+(p-1)(h^2+h')v\big)<0\Big\},\endaligned$$
In Theorem \ref{CNth1} below, for each $f\in\scr{F}_I$,
$\inf_{x\in(0,D)}I(f)(x)^{-1}$ produces a lower bound of
$\lambda_p$. So the part having ``$\sup\inf$'' in each of the
formulas is used for the lower estimates of $\lambda_p$. Dually, the
part having ``$\inf\sup$'' is used for the upper estimates. These
formulas deduce the basic estimates in Theorem \ref{CNth2} and the
approximating procedure in Theorem \ref{CNth3}.
\begin{thm}\label{CNth1}$(\text{\rm Variational formulas})$  For $p>1$, we have
\begin{itemize}\setlength{\itemsep}{-0.8ex}
\item[$(1)$] single integral forms:
$$\aligned
\inf_{f\in\mathscr{{\widetilde F}}_I}\sup_{x\in(0,D)}I(f)(x)^{-1}
=\lambda_p=\sup_{f\in\mathscr F_I}\,\inf_{x\in(0,D)}I(f)(x)^{-1},
\endaligned$$
\item[$(2)$]double integral forms:
$$
\aligned
\inf_{f\in\mathscr{{\widetilde F}}_{I\!I}}\,\sup_{x\in
{\supp}(f)}I\!I(f)(x)^{-1}=\lambda_p=\sup_{f\in\mathscr{F}_{I\!I}}\,
\inf_{x\in(0, D)}I\!I(f)(x)^{-1}.
\endaligned$$
\end{itemize}
Moreover, if $u$ and $v'$ are continuous,
 then we have additionally
\begin{itemize}\setlength{\itemsep}{-0.8ex}
\item[$(3)$]  differential forms:
$$
\aligned
 \inf_{h\in\widetilde{\mathscr{H}}}\,\sup_{x\in(0,
D)}R(h)(x)=\lambda_p=\sup_{h\in\mathscr{H}}\,\inf_{x\in(0,
D)}R(h)(x).
\endaligned$$
\end{itemize}Furthermore,
the supremum on the right-hand side of the above three formulas can
be attained.
\end{thm}

The following proposition adds some additional sets of functions for
operators $I$ and $I\!I$. It then provides alternative descriptions
of the lower and upper estimates of $\lambda_p$.
\begin{prop}\label{CNprop}For $p>1$, we have
$$\aligned
\lambda_p=\sup_{f\in\mathscr{F}_{I}}\inf_{x\in(0, D)}I\!I(f)(x)^{-1};
\endaligned
$$
$$\aligned
\!\!\lambda_p&\!=\!\!\!\!\!\inf_{f\in
\widetilde{\mathscr{F}}_{I\!I}\cup\widetilde{\mathscr{F}}_{I\!I}'}\sup_{x\in\,{\supp}(f)}\!\!\!\!\!
I\!I(f)(x)^{-1}
\!\!=\!\!\!\inf_{f\in
\widetilde{\mathscr{F}}_{I}}\sup_{x\in\,{\supp}(f)}\!\!\!\!\!I\!I(f)(x)^{-1}
\!\!=\!\!\!\inf_{f\in\widetilde{\mathscr{F}}_{I}'}\sup_{x\in(0,
D)}\!\!\!I(f)(x)^{-1}\!\!,\!\!\!
\endaligned
$$
where
$$\widetilde{\mathscr{F}}_{I\!I}'\!=\!\big\{f:  f\in\mathscr{ C}[0, D] \text{ and }
fI\!I(f)\in L^p(\mu) \big\},$$
$$\aligned&\widetilde{\mathscr{F}}_{I}'\!=\!\big\{f:\exists\,x_{0}\in(0, D), f\!=\!f
\mathbbold{1}_{[0, x_{0})}\in\mathscr{C}[0, x_0], f'|_{(0, x_{0})}\!<\!0,\\
&\hskip1.7cm \text{and } v^{p^*-1}f'\!\in\! \mathscr{C}(0, x_{0})
\big\}.\endaligned$$
\end{prop}

Define $k(p)=p{p^*}^{p-1}$ for $p>1$ and
$$\sigma_p=\sup_{x\in(0,D)}\mu(0,x)\, \hat{\nu}(x, D)^{p-1}.$$
As applications of the variational formulas in Theorem \ref{CNth1}
(1), we have the following basic estimates
known in \cite{BA1990}.
\begin{thm}\label{CNth2}$(\text{\rm Criterion and basic estimates})$
For $p>1$, the eigenvalue $\lambda_p>0$ if and only if
$\sigma_p<\infty$. Moreover, the following basic estimates hold:
 \begin{equation*}
 (k(p)\sigma_p)^{-1}\leqslant\lambda_{p}\leqslant\sigma_p^{-1},
 \end{equation*}
 In particular, we have $\lambda_p=0$ if ${\hat\nu}(0, D)=\infty$ and $\lambda_p>0$ if
 $$\int_0^{D}\mu(0, s)^{p^*-1}{\hat v}(s)\d s<\infty.$$
 \end{thm}

 The approximating procedure below is  an application of variational formulas in Theorem \ref{CNth1}\,(2).
 The main idea is an iteration, its first step produces Corollary \ref{CNcor1} below. Noticing that
 $\lambda_p$ is trivial once $\sigma_p=\infty$, we may assume that $\sigma_p<\infty$ for further
 study on the estimates of $\lambda_p$.
 \begin{thm}\label{CNth3}
 $(\text{\rm Approximating procedure})$
Assume that $\sigma_p<\infty$.
\begin{itemize}\setlength{\itemsep}{-0.8ex}
\item[$(1)$] Let $f_1=\hat{\nu}(\cdot, D)^{1/p^*}$,
$f_{n+1}=f_{n}I\!I(f_{n})^{p^*-1}$ and $\delta_{n}=\sup_{x\in(0,D)}I\!I(f_{n})(x)$ for  $n\geqslant 1$.
 Then $\delta_{n}$ is decreasing and $$\lambda_{p}\geqslant\delta_{n}^{-1}\geqslant(k(p)\sigma_p)^{-1}.$$
\item[$(2)$]  For fixed $x_{0}$, $x_{1}\in(0,D)$ with $x_{0}<x_{1}$, define
 $$f_1^{x_{0},x_{1}}=\hat\nu(\cdot\vee x_{0}, x_1)\mathbbold{1}_{ [ 0, x_{1})},\qquad
 f_{n}^{x_{0},x_{1}}=f_{n-1}^{x_{0},x_{1}} I\!I\big(f_{n-1}^{x_{0},x_{1}}\big)^{p^*-1}
 \mathbbold{1}_{[0,x_{1})},$$
and
$$\delta_{n}'=\sup_{x_0,x_1:\,x_{0}<x_{1}}\inf_{x< x_{1}}
I\!I(f_{n}^{x_{0},x_{1}})(x),\quad \text{for}\ n\geqslant 1.$$
 Then $\delta_{n}'$ is increasing and
 $$\sigma_p^{-1}\geqslant{\delta_{n}'}^{-1}\geqslant\lambda_{p}.$$
Next, define $$\bar{\delta}_n=\sup_{x_{0}<x_{1}}\frac{\|f_n^{x_{0},x_{1}}\|_p^p}{D_{p}(f_{n}^{x_{0},x_{1}})},
\qquad n\geqslant 1.$$
 Then $\bar{\delta}_n^{-1}\geqslant\lambda_{p}$ and $\bar\delta_{n+1}\geqslant \delta_n'$ for $n\geqslant1$.
\end{itemize}
\end{thm}

The following Corollary \ref{CNcor1} can be obtained directly from
Theorem \ref{CNth3}. It provides us some improved and explicit
estimates of the eigenvalue (see Example \ref{CNEX1} below).
\begin{cor}\label{CNcor1}$(\text{\rm Improved estimates})$Assume that  $\sigma_p<\infty$. Then
$$\sigma_p^{-1}\geqslant{\delta_1'}^{-1}\geqslant\lambda_p\geqslant\delta_1^{-1}\geqslant(k(p)\sigma_p)^{-1}
,$$
where
$$\delta_1=\sup_{x\in(0, D)}\bigg[\frac{1}{\hat{\nu}(x, D)^{1/p^*}}\int_x^D\hat{v}(s)
\bigg(\int_0^s\hat{\nu}(t, D)^{(p-1)/p^*}\mu(\d t)\bigg)^{p^*-1}\d s\bigg]^{p-1};$$
$$\aligned
\delta_1'&=\sup_{x\in (0, D)}\frac{1}{\hat{\nu}(x, D)^{p-1}}\bigg[\int_{x}^{D}\hat{v}(s)
\bigg(\int_0^s\hat{\nu}(t\vee x, D)^{p-1}\mu(\d t)\bigg)^{p^*-1}\d s\bigg]^{p-1}.\endaligned$$
Moreover,
$$\bar\delta_1=\sup_{x\in(0, D)}\bigg[\mu(0,x)\hat{\nu}(x, D)^{p-1}+
\frac{1}{\hat{\nu}(x, D)}\int_{x}^D\hat{\nu}(t, D)^p\mu(\d t)\bigg]\in[\sigma_p,  p \sigma_p],$$
 and $\bar\delta_1\leqslant{\delta_1'}$ for $1<p\leqslant2$, $\bar\delta_1\geqslant{\delta_1'}$
for $p\geqslant2$.
\end{cor}

When $p=2$, the assertion that $\bar\delta_1=\delta_1'$ was proved
in \rf{CWZ2013}{Theorem 3}. To illustrate the results above, we present an
example as follows.

\begin{exm}\label{CNEX1}{\rm
Let $\d\mu=\d\nu=\d x$ on $(0, 1)$. In the ND-case,
the eigenvalue $\lambda_p$ is
\begin{equation}
\lambda_p^{1/p}=\frac{\pi (p-1)^{1/p}}{p}\sin^{-1}\frac \pi p.\end{equation}
For the basic estimates, we have
$$\sigma_p^{1/p}=\bigg(\frac 1 p\bigg)^{1/p}
\bigg(\frac 1 {p^*}\bigg)^{1/p^*}.$$
Furthermore, we have
\begin{gather}{\bar\delta}_1^{1/p}=p^{1/p - 2} \big(p^2-1\big)^{1-1/p},\nonumber\\
\delta_1^{1/p}=\frac{1}{(p+1/p-1)^{1/p}}\bigg\{\sup_{x\in (0, 1)}\frac{1}{(1-x)^{1/p^*}}
\int_0^{1-x}\big(1-z^{p+1/p-1}\big)^{p^*-1}\d
z\bigg\}^{1/p^*}.\nonumber\end{gather}
\begin{figure}[h]
 \begin{center}{\includegraphics[width=11.0cm,height=6.45cm]{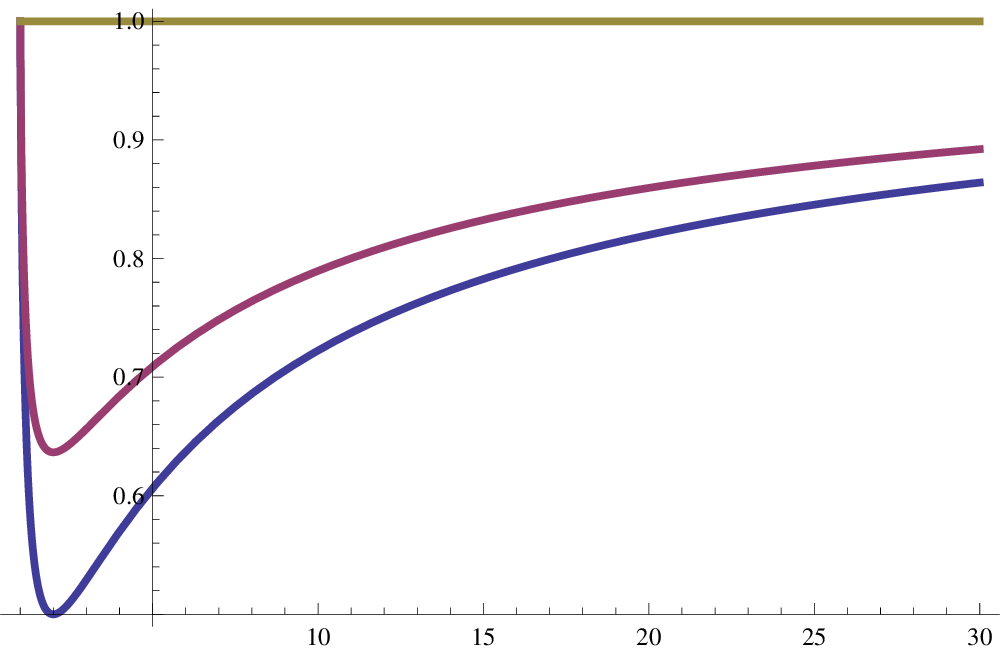}\newline
{\bf Figure 1}\quad The middle curve is the exact value of
$\lambda_p^{1/p}$. The top straight line and the bottom curve are
the basic estimates of $\lambda_p^{1/p}$.}\end{center}
\end{figure}\vspace{-0.8truecm}
\begin{figure}[h]
\begin{center}{\includegraphics[width=11.0cm,height=6.45cm]{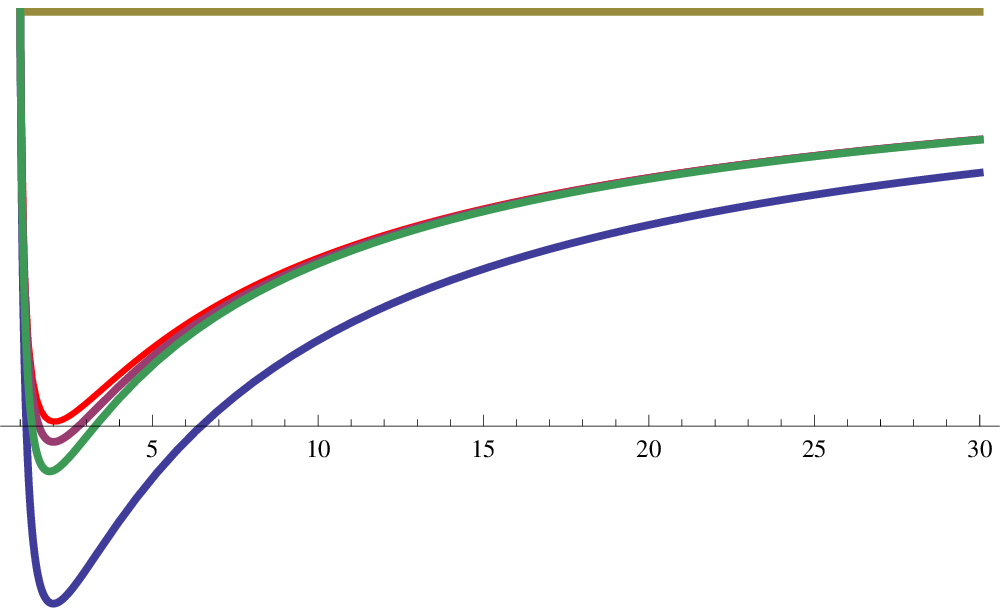}\newline
{\bf Figure 2}\quad The improved bounds ${\delta}_1^{1/p}$ and
${\bar\delta}_1^{1/p}$ are added to Figure 1.}\end{center}
\end{figure}
\noindent The exact value $\lambda_p^{1/p}$ and its basic estimates
are shown in Figure 1. Then, the improved upper bound
${\delta}_1^{1/p}$ and lower one ${\bar\delta}_1^{1/p}$ are added to
Figure 1, as shown in Figure 2. It is quite surprising and
unexpected that both of ${\delta}_1^{1/p}$ and
${\bar\delta}_1^{1/p}$ are almost overlapped with the exact value
$\lambda_p^{1/p}$ except in a small neighborhood of $p=2$, where
$\delta_1^{1/p}$ is a little bigger and ${\bar\delta}_1^{1/p}$ is a
little smaller than $\lambda_p^{1/p}$. Here ${\delta_1'}^{1/p}$ is
ignored since it improves ${\bar\delta}_1^{1/p}$ only a little bit
for $p\in (1, 2)$. }
\end{exm}
\noindent\\[4mm]

\setcounter{section}{3}
\setcounter{thm}{0}

\noindent{\bbb 3\quad Proofs of the main results}\\[0.1cm]
\noindent Some preparations for the proofs are collected in Subsection
3.1. They may not be used completely in the proofs but are
helpful to understand the idea in this paper and may be useful in
other cases. The proofs of the main results are presented in
Subsection 3.2. For simplicity, we let $\uparrow$ (resp.
$\upuparrows$, $\downarrow$, $\downdownarrows$) denote
increasing (resp. strictly increasing, decreasing, strictly
decreasing) throughout this paper.

\vspace{2mm}
\noindent{\bbb3.1\quad Preparations}

\vspace{1mm}\noindent
The next lemma is taken from \rf{CDF2004}{Theorem 1.1 on page 170}(see \cite{W1947} for its original idea).
Combining with the following Remark \ref{CND-rem1}, Lemmas \ref{CNL0} and \ref{CNL01},
it guarantees the existence of the solution $(\lambda_p, g)$ to the  eigenvalue problem.
\begin{lem}$($\text{\rm Existence and Uniqueness}$)$
 \begin{itemize}\setlength{\itemsep}{-0.8ex}
 \item[$(1)$]Suppose that  $u$
and $v$ are locally integrable on $[0, D]\subseteq\mathbb{R}$ $($or $[0, D)\subseteq\mathbb{R}$ provided $D=\infty$$)$ and $v>0$. Given constants $A$ and $B$, for each fixed $\lambda$, there is
uniquely a solution $g$ such that $g(0)=A$, $g'(0)=B$ and  the eigenequation \eqref{CND-eigen} holds almost everywhere. Moreover, $v^{p*-1}g'$ is absolutely continuous.
\item[$(2)$]Suppose additionally $u$ and $v$ are continuous. Then $g\in\mathscr{C}^2[0, D]$ and the eigenequation holds everywhere on $[0, D]$.
 \end{itemize} \label{Cexists}\end{lem}
 \vspace{-1ex}
 If the eigenequation \eqref{CND-eigen} holds (almost) everywhere for $(\lambda_p, g)$, then $g$ is called an (a.e.) eigenfunction of $\lambda_p$.
\begin{rem}\label{CND-rem1}
$(1)$ One may also refer to  {\rm\rf{DK2005}{Lemma 2.1}} for the existence of solution to eigenvalue problem with ND boundary conditions provided $D<\infty$. When $D=\infty$,  the Dirichlet boundary at  $D$ means $g(D)=0$, which is proved by Proposition $\ref{CNp1}$ below.

$(2)$ By {\rm\rf{BA1990}{Theorem 4.1, Theorem 4.7}}, we see that the eigenequation in \eqref{CND-eigen} has  solutions if and if only the following equation has solutions:
\begin{equation*}\big(|g'|^{p-2}g'\big)'(x)=-\lambda \tilde u(x)|g|^{p-2}g(x)\end{equation*}
where $\tilde u$ is related to $v$ and $u$ in the eigenequation. Hence the weight function $v$ in the eigenquation is not a sensitive or key quantity to the existence of solution to the eigenequation and can be seen as a constant.
\end{rem}
Define $\scr{A}_K[0, D]=\{f: f\in\scr{A}[0, D], f \text{ has compact
support}\}$ and
\begin{equation}\lambda_{*, p}=\inf\{D_p(f): f\in\scr{A}_K[0, D], \|f\|_p=1,
\text{ and }f(D)=0 \text{ if } D<\infty\},
\label{CNDf2}\end{equation}
\begin{equation}\tilde\lambda_p=\inf\{D_p(f): v^{p*-1}f'\in\scr{C}(0, D), f\in\scr{C}[0,D],
\|f\|_p=1, f(D)=0 \},
\label{CNDf3}\end{equation} where $\|\cdot\|_p$ means the norm in $L^p(\mu)$ space. The following quantities are also useful
for us. Set $\alpha\in(0, D)$ and define
$$\aligned&\lambda_{*,p}^{(0, \alpha)}=\inf\big\{D_p(f): \mu\big(|f|^p\big)=1,
f\in\mathscr{A}[0, \alpha]  \text{ and }f|_{[\alpha, D]}=0\big\}.\endaligned$$
$$\lambda_p^{(0, \alpha)}
=\inf\!\big\{D_p\big(f\big)\!: f\in\mathscr{C}[0,\alpha], v^{p^*-1}f'\in\mathscr{C}(0,\alpha), \mu\big(|f|^p\big)\!=\!1, f|_{[\alpha,D]}\!=\!0\big\}.
$$
The following three Lemmas describe in a refined way the first
eigenvalue and lead to, step by step, the conclusion that
$$\tilde\lambda_p=\lambda_p=\lambda_{*,p}=\tilde\lambda_{*,p}.$$

\begin{lem}\label{CNL0} We have $\lambda_p=\lambda_{*, p}$.
\end{lem}
\noindent{\it Proof}\quad It is obvious that
$\lambda_p\geqslant\lambda_{*, p}$.  Next, let $g$ be the a.e.
eigenfunction of $\lambda_{*,p}$. Then $g\in\mathscr{C}[0,D]$ and $v^{p*-1}g'\in\scr{C}(0, D)$ by
Lemma \ref{Cexists}. Since $L_pg=-\lambda_{*,p}|g|^{p-2}g$,
by the arguments after formula \eqref{CNDf1}, we have
$$-\big(vg|g'|^{p-2}g'\big)\big|_{0}^{D}+D_p(g)=\lambda_{*,p}\|g\|_p^{p}.$$
Since $g'(0)=0$ and $(gg')(D)\leqslant 0$, we have
$\lambda_{*,p}\geqslant D_p(g)/\|g\|_p^p$. Because
$g\in\mathscr{C}_K[0,D]$, it is clear that
$D_p(g)/\|g\|_p^p\geqslant\lambda_p$. We have thus obtained that
$$\lambda_p\leqslant\lambda_{*,p}\leqslant\lambda_p,$$ and so
$\lambda_p=\lambda_{*,p}$. There is a small gap in the proof above
since in the case of $D=\infty$, the a.e. eigenfunction $g$ may not
belong to $L^p(\mu)$ and we have not yet proved that
$(gg')(D)\leqslant 0$. However, one may avoid this by a standard
approximating procedure,
using $[0,\alpha_{n}]$ instead of $[0,D)$ with $\alpha_{n}\uparrow D$ provided $D=\infty$:
$$\aligned
\lim_{n\to\infty}\lambda_p^{(0, \alpha_n)}=\lim_{n\to\infty}\inf\big\{&D_p(f):\mu\big(|f|^p\big)\!\!=\!\!1,
f\in\!{\scr C}[0,\alpha_n], v^{p^*-1}f'\in\scr{C}(0, \alpha_n),\\
& f|_{[\alpha_n, D]}\!=\!0\big\}\!=\!\!\lambda_p.\endaligned$$
Similarly, $\lambda_{*,p}^{(0,\alpha_n)}\to\lambda_{*,p}$ as $n\to \infty$.$\qquad\Box$
\begin{lem}\label{CNL01} For $\tilde\lambda_{*,p}$ defined in \eqref{CNDf1}, we have
$\tilde\lambda_{*,p}=\lambda_{*,p}$. Furthermore,
$\tilde\lambda_p=\lambda_p=\lambda_{*,p}=\tilde\lambda_{*,p}$.
\end{lem}
{\noindent\it Proof}\quad On one hand, by definition, if $\beta_{n+1}\geqslant \beta_n$,
then $\lambda_{*,p}^{(0, \beta_n)}\geqslant \lambda_{*,p}^{(0, \beta_{n+1})}$. We have thus obtained
$$\lim_{n\to\infty} \lambda_{*,p}^{(0, \beta_n)}\geqslant \lambda_{*,p}^{(0, D)}
= {\tilde\lambda}_{\ast,p}.$$
On the other hand, by definition of $\tilde\lambda_{*,p}$, for any
fixed $\varepsilon>0$, there exists $f$ satisfying $\|f\|_p=1$,
$f(D)=0$, and $D_p(f)\leqslant
{\tilde\lambda}_{\ast,p}+\varepsilon$. Let $\beta_n\uparrow D$ and
$f_n=\big(f-f(\beta_n)\big)\mathbbold{1}_{[0, \beta_n)}$. Then
$D_p(f_n)\uparrow D_p(f)$ as $n\uparrow\infty$. Choose subsequence
$\{n_m\}_{m\geqslant 1}$ if necessary such that
$$\varlimsup_{n\to\infty}\frac{D_p(f_n)}{\|f_n\|_p^p}
=\lim_{m\to\infty}\frac{D_p(f_{n_m})}{\|f_{n_m}\|_p^p}.$$
By Fatou's lemma and the fact that $f(D)=0$, we have
$$\varliminf_{m\to\infty}{\|f_{n_m}\|_p^p}\geqslant \Big\|\varliminf_{m\to\infty}f_{n_m}\Big\|_p^p
=\|f\|_p^p=1.$$
Therefore, we obtain
$$\aligned\varlimsup_{n\to\infty}\lambda_{*,p}^{(0, \beta_n)}
&\leqslant \varlimsup_{n\to\infty}\frac{D_p(f_n)}{\|f_n\|_p^p}
=\lim_{m\to\infty}\frac{D_p(f_{n_m})}{\|f_{n_m}\|_p^p}
\leqslant \frac{\lim_{m\to\infty}D_p(f_{n_m})}{\varliminf_{m\to\infty}\|f_{n_m}\|_p^p}\!\! \leqslant D_p(f)\\
&\leqslant \tilde\lambda_{*,p}+\varepsilon.
\endaligned$$
Since $\lim_{n\to\infty}\lambda_{*,p}^{(0, \beta_n)}=\lambda_{*,p}$,
we get $\tilde\lambda_{*,p}=\lambda_{*,p}$.
 Moreover,   $$\tilde\lambda_p\geqslant\tilde\lambda_{*,p}=\lambda_{*,p}=\lambda_p\geqslant\tilde\lambda_p$$
 and the required assertion holds. $\qquad\Box$

The following lemma, which serves for Lemma \ref{CNJM20122}, presents us that $\{\lambda_{*,p}^{(0, \alpha)}\}$ is strictly decreasing with respect to $\alpha$.
\begin{lem}\label{C-NL1} For $\alpha, \beta\in(0, D)$ with $\alpha<\beta$,
we have $\lambda_{*, p}^{(0, \alpha)}>\lambda_{*,p}^{(0, \beta)}$. Furthermore,
$\lambda_{*, p}^{(0, \beta_n)}\downdownarrows \lambda_{*, p}$ as $\beta_n\upuparrows D$.
\end{lem}
\noindent{\it Proof}\quad\rm  Let $g\,(\ne 0)$ be an a.e. eigenfunction of $\lambda_{*, p}^{(0, \alpha)}$.
Then $g'(0)=0$, $g(\alpha)=0$, and $L_pg=-\lambda_{*, p}^{(0, \alpha)} |g|^{p-2}g$  on $(0, \alpha)$. Moreover,
$$\lambda_{*, p}^{(0, \alpha)}=\frac{D_p^{0, \alpha}(g)}{\|g\|_{L^p(0, \alpha;\, \mu)}^p},
\qquad D_p^{\alpha, \beta}(f)=\int_{\alpha}^\beta  |f'|^p \d\nu$$
(see arguments after formula \eqref{CNDf1}).
By the proof of Lemma \ref{CNL0}, the proof of the first assertion  will be done once we choose a
function $\tilde g\in\mathscr{A}[0,\beta]$
such that $\tilde g'(0)=0$, $\tilde g(\beta)=0$, and
\begin{equation}\label{CNDff3}\frac{D_p^{0, \alpha}(g)}{\|g\|_{L^p(0, \alpha;\, \mu)}^p}
>\frac{D_p^{0, \beta}(\tilde g)}{\|\tilde g\|_{L^p(0, \beta;\, \mu)}^p}\quad
\Big(\geqslant \lambda_{*, p}^{(0, \beta)}\Big).\end{equation}
To do so, without loss of generality, assume that $g|_{(0, \alpha)}>0$ (see \rf{P2006}{Lemma 2.4}).
Then the required assertion follows
for
$$ {\tilde g}(x)=
(g+\varepsilon)\mathbbold{1}_{[0, \alpha)}(x)
+\frac{\varepsilon(\beta-x)}{(\beta-\alpha)}\mathbbold{1}_{[\alpha, \beta]}(x),\qquad x\in[0, \beta],
$$ once $\varepsilon$ is sufficiently small.
Actually, by simple calculation, we have
$$\aligned &D_p^{0, \beta}(\tilde g)=D_p^{0, \alpha}(g)
+\frac{\varepsilon^{p}}{(\beta-\alpha)^p}\nu(\alpha, \beta),\endaligned$$
$$\aligned\|\tilde g\|^p_{L^p(0,\beta;\mu)}=\|g\|^p_{L^p(0, \alpha; \mu)}
+\int_{0}^{\alpha}\big(|g+\varepsilon|^p-|g|^p\big)\d\mu
+ \int_{\alpha}^{\beta}\frac{\varepsilon^p(\beta-x)^p}{(\beta-\alpha)^p}\mu(\d x).\endaligned$$
Since $$\lambda_{*, p}^{(0, \alpha)}= D_p^{0,
\alpha}(g)\big/\|g\|_{L^p(0, \alpha;\mu)}^p,$$ inequality
\eqref{CNDff3} holds if and only if
$$\aligned&\frac{\varepsilon^p\nu(\alpha, \beta)}{(\beta-\alpha)^p}<\bigg(\int_{0}^{\alpha}
\big(|g+\varepsilon|^p-|g|^p\big)\d\mu+\frac{\varepsilon^p}{(\beta-\alpha)^p}
\int_{\alpha}^{\beta}(\beta-x)^p\mu(\d x)\bigg)\lambda_{*,p}^{(0, \alpha)}.\endaligned$$
It suffices to show that
$$\frac{\varepsilon^{p-1}}{(\beta-\alpha)^p}\nu(\alpha, \beta)<\lambda_{*, p}^{(0, \alpha)}
\bigg(\int_{0}^{\alpha}\frac{|g(x)+\varepsilon|^p-|g(x)|^p}{\varepsilon}\mu(\d x)\bigg).$$
By letting $\varepsilon\to 0$, the right-hand side is equal to
$$\lambda_{*, p}^{(0, \alpha)}\int_{0}^{\alpha}pg^{p-1}\d\mu,$$
which is positive. So the required inequality is obvious for sufficiently small
$\varepsilon$ and the first assertion holds.
The second assertion  was proved at the end of the proofs of Lemma \ref{CNL01}. $\qquad\Box$

The following Lemma is about the  eigenfunction of $\lambda_p$,
which is the basis of the test functions used for the corresponding
operators.
\begin{lem}\label{CNJM20122}
Let $g$ be the first eigenfunction of eigenvalue problem \eqref{CND-eigen}. Then both
$g$ and $g'$ do not change sign. Moreover, if $g>0$, then $g'<0$.
\end{lem}
{\noindent\it Proof}\quad If there exists $\alpha\in(0, D)$ such
that $g(\alpha)=0$, then $\lambda_{*,p}^{(0,
\alpha)}\leqslant\lambda_{*, p}$ by the minimum property of
$\lambda_{*,p}^{(0, \alpha)}$. However, by Lemma \ref{C-NL1}, we get
$\lambda_{*,p}^{(0, \alpha)}\downdownarrows \lambda_{*,p}$ as
$\alpha\upuparrows D$. This is a contradiction. So $g$ does not
change its sign. Next, consider $g'$. By \rf{P2006}{Lemma 2.3}, if
there exists $x\in(0, D)$ such that $g'(x)=0$, then $\exists
x_0\in(0, x)$ such that $g(x_0)=0$, which is impossible by the
strictly decreasing property of $\lambda_{*,p}^{(0, \alpha)}$ with
respect to $\alpha$. So the assertion holds.$\qquad\Box$

Before moving on, we introduce a general equation, non-linear `Poisson equation' as follows:
\begin{equation}\label{C-Poisson}L_pg(x)=-u(x)|f|^{p-2}f(x),\qquad x\in(0, D).\end{equation}
Integration by parts yields that for $x, y\in(0, D)$ with $x<y$,
\begin{equation}
 \label{C-N-D}v(x)|g'|^{p-2}g'(x)- v(y)|g'|^{p-2}g'(y)=\int_x^y|f|^{p-2}f\d\mu.\end{equation}
 By replacing $f$ with $\lambda^{p^*-1}g$,  it is not hard to understand where the operator $I$  comes from.
Moreover, if $g$ is positive and decreasing, $g'(0)=0$, then
\begin{equation}
\label{C-Nf1}g(y)-g(D)
=\int_y^D\hat{v}(x)\bigg(\int_0^x|f|^{p-2}f\d\mu\bigg)^{p^*-1}\d x,\qquad y\in(0, D).
\end{equation}
By replacing $f$ with $\lambda^{p^*-1}g$,  it is easy to see where the
operator $I\!I$ comes from, provided $g(D)=0$ (which is affirmative
by Proposition \ref{CNp1} below). Finally, assume that $(\lambda_p,
g)$ is a solution to \eqref{CND-eigen}. Then
$\lambda_p=-{L_pg}/{\big(|g|^{p-2}g u\big)}$. Hence, by letting
$h=g'/g$, we deduce the operator $R$ from the eigenequation.

\vspace{2mm}
\noindent{\bbb3.2\quad
Proof of the main results}

\vspace{1mm}\noindent
{\it Proof of Theorem $\ref{CNth1}$ and Proposition $\ref{CNprop}$}\quad We adopt the circle
arguments below to prove the lower estimates:
$$\begin{aligned}\lambda_{p}\geqslant\tilde\lambda_p&\geqslant\sup_{f\in\mathscr{F}_{I\!I}}
\inf_{x\in(0, D)}I\!I(f)(x)^{-1}
=\sup_{f\in\mathscr{F}_{I}}\inf_{x\in(0, D)}I\!I(f)(x)^{-1}\\&\hskip1cm=
\sup_{f\in\mathscr{F}_{I}}\inf_{x\in(0, D)}I(f)(x)^{-1}\\&\geqslant
\sup_{h\in\mathscr{H}}\inf_{x\in(0, D)}R(h)(x)\\
&\geqslant\lambda_p. \end{aligned}$$

Step 1\quad Prove that
$\lambda_{p}\geqslant\tilde\lambda_p\geqslant\sup_{f\in\mathscr{F}_{I\!I}}\inf_{x\in(0,
D)}I\!I(f)(x)^{-1}.$

It suffices to show the second inequality. For each fixed  $h>0$ and
$g\in\mathscr{C}[0, D]$ with $\|g\|_{p}=1$, $g(D)=0$ and $v^{p^*-1}g'\in \mathscr{C}(0, D)$,  we have
$$\begin{aligned}
  \int_{0}^{D}|g^p|\d\mu&=\int_{0}^{D}\bigg|\int_{x}^{D}g'(t)\bigg(\frac{v(t)}{h(t)}
  \bigg)^{1/p}\bigg(\frac{h(t)}{v(t)}\bigg)^{1/p}\d t\bigg|^{p}\mu(\d x)\\
&\leqslant
\int_{0}^{D}\int_{x}^{D}\frac{v(t)}{h(t)}\big|g'(t)\big|^p\d t\bigg[\int_{x}^{D}\bigg(\frac{h(s)}{v(s)}
\bigg)^{p^*-1}\d s\bigg]^{p-1}\mu(\d x)\\
&\hskip3cm\quad(\text{by H\"{o}lder's inequality})\\
&=
\int_{0}^{D}\frac{v(t)}{h(t)}\big|g'(t)\big|^p\d t\int_{0}^{t}\bigg[\int_{x}^{D}\bigg(\frac{h(s)}{v(s)}
\bigg)^{p^*-1}\d s\bigg]^{p-1}\mu(\d x)\\&\hskip3cm\quad (\text{by Fubini's Theorem})\\
&\leqslant D_p(g)\sup_{t\in(0, D)}H(t),
 \end{aligned}$$
 where
 $$H(t)=\frac{1}{h(t)}\int_{0}^{t}\bigg[\int_{x}^{D}\bigg(\frac{h(s)}{v(s)}\bigg)^{p^*-1}\d s\bigg]^{p-1}
 \mu(\d x).$$
For  $f\in \mathscr{F}_{I\!I}
$ with $\sup_{x\in(0, D)}I\!I(f)(x)<\infty$, let
$$h(t)=\int_{0}^{t}f^{p-1}(s)u(s)\d s.$$ Then
 $h'=f^{p-1}u$. By  Cauchy's mean-value theorem, we have
  $$\aligned
\sup_{x\in(0, D)}\!H(x)&\leqslant\sup_{x\in(0, D)}\!I\!I(f)(x). \endaligned$$
Thus\ $\lambda_{p}\geqslant\inf_{x\in(0, D)}I\!I(f)(x)^{-1}$. The assertion then follows by making
the supremum with respect to $f\in\mathscr{F}_{I\!I}$.

Step 2\quad Prove that
 $$\begin{aligned}
\sup_{f\in\mathscr{F}_{I\!I}}\inf_{x\in(0, D)}I\!I(f)(x)^{-1}=
\sup_{f\in\mathscr{F}_{I}}\inf_{x\in(0, D)}I\!I(f)(x)^{-1}=
\sup_{f\in\mathscr{F}_{I}}\inf_{x\in(0, D)}I(f)(x)^{-1}.
\end{aligned}$$

 (a)\quad We prove the part `$\geqslant$'.  Since
$\mathscr{F}_{I}\subset\mathscr{F}_{I\!I}$, it suffices to show that
$$\sup_{f\in\mathscr{F}_{I}}\inf_{x\in(0,
D)}I\!I(f)(x)^{-1}\geqslant \sup_{f\in\mathscr{F}_{I}}\inf_{x\in(0,
D)}I(f)(x)^{-1}$$
for $f\in\mathscr{F}_{I}$ with $\sup_{x\in(0,
D)}I(f)<\infty$. Since $f(D)\geqslant0$, by replacing $f$ in the
denominator of $I\!I(f)$  with $-\int_{\cdot}^D f'(s)\d s$ and using
 Cauchy's mean-value theorem, we have
$$\aligned
\sup_{x\in (0, D)}I\!I(f)(x)&\leqslant\sup_{x\in(0, D)}I(f)(x)<\infty.
\endaligned$$
 So the assertion holds by making the supremum with respect to $f\in\mathscr{F}_I$.

(b)\quad To prove the equality, it suffices  to show that
 $$ \sup_{f\in\mathscr{F}_{I}}\inf_{x\in(0,D)}I(f)(x)^{-1}\geqslant
\sup_{f\in\mathscr{F}_{I\!I}}\inf_{x\in(0, D)}I\!I(f)(x)^{-1}.$$
For  $ f\in\mathscr{F}_{I\!I}$,
without loss of generality, assume that
$\inf_{x\in(0, D)}I\!I(f)(x)^{-1}>0$.
Let
$g=f[I\!I(f)]^{p^*-1}$. Then $g\in\mathscr{F}_{I}$. Moreover,
$$v(x)\big(-g'(x)\big)^{p-1}
=\int_{0}^{x}f^{p-1}\d\mu\geqslant
\int_0^xg^{p-1}\d\mu\inf_{t\in(0,x)}\frac{f^{p-1}(t)}{g^{p-1}(t)},$$
i.e.,
$$I(g)(x)^{-1}\geqslant \inf_{x\in(0,D)}I\!I(f)(x)^{-1}.$$
Hence,
\begin{equation*}
\sup_{f\in\mathscr{F}_{I}}\inf_{x\in(0, D)}I(f)(x)^{-1}\geqslant
\inf_{x\in(0, D)}I(g)(x)^{-1}\geqslant \inf_{x\in(0,
D)}I\!I(f)(x)^{-1}
\end{equation*}
and the assertion holds since
$f\in\mathscr{F}_{I\!I}$ is arbitrary.

Then there is another method to prove the equality: prove that
$$\sup_{f\in\mathscr{F}_I}\inf_{x\in(0, D)}I(f)(x)^{-1}\geqslant\lambda_p.$$
Let $g$ be an a.e. eigenfunction corresponding to  $\lambda_p$.
Then $g$ is positive and strictly decreasing. It is easy to check that $g\in\mathscr{F}_{I}$.
By \eqref{C-N-D}, we have $$\lambda_p=\inf_{x\in(0,D)}I(g)(x)^{-1}
\leqslant\sup_{f\in\mathscr{F}_I}\inf_{x\in(0,D)}I(f)(x)^{-1}.$$

Step 3\quad When $u$ and $v'$ are continuous, we prove that
$$\sup_{f\in\mathscr{F}_{I\!I}}\inf_{x\in(0,
D)}I\!I(f)(x)^{-1}\geqslant \sup_{h\in\mathscr{H}}\inf_{x\in(0,
D)}R(h)(x). $$

First, we change the form of $R(h)$. Let $g$ with $g(D)=0$ be a positive function on $[0, D)$ such that
$h=g'/g$ (see the arguments after Lemma \ref{CNJM20122}).
Then
$$\aligned R(h)=-u^{-1}\big\{|h|^{p-2}\big[v' h+(p-1)v(h^2+h')\big]\big\}
=-\frac{1}{u g^{p-1}}L_pg.\endaligned$$
Now, we turn to our main text. It suffices to show that
$$\sup_{f\in\mathscr{F}_{I\!I}}\inf_{x\in(0,D)}I\!I(f)(x)^{-1}\geqslant \inf_{x\in(0, D)}R(h)(x)
\quad\text{ for every } h\in\mathscr{H}.$$
 Without loss of generality, assume that $\inf_{x\in(0, D)}R(h)(x)>0$,
which implies $R(h)>0$ on $(0,D) $. Let $f=g(R(h))^{p^*-1}$ ($g$ is
the function just specified). Since $u,$ $v'$ are continuous, we
have $f\in\mathscr{F}_{I\!I}$ and
\begin{equation*}
u(x) f^{p-1}(x)=-L_{p}g(x),\qquad  x\in (0, D).
\end{equation*}
Moreover, by \eqref{C-Nf1}, we have
$$\aligned
 g(y)-g(D)&=\int_{y}^{D}\hat{v}(x)\bigg(\int_{0}^{x}f^{p-1}\d\mu\bigg)^{p^*-1}\d x.
\endaligned$$
So
${g^{p-1}}/{f^{p-1}}\geqslant I\!I(f)$ on $(0, D)$ and \begin{equation*}
\inf_{(0, D)}R(h)=\inf_{(0, D)}f^{p-1}/{g^{p-1}}\leqslant\inf_{(0,
D)}I\!I(f)^{-1}\leqslant\sup_{f\in\mathscr{F}_{I\!I}}\inf_{x\in(0,
D)}I\!I(f)(x)^{-1}.
\end{equation*}
 Hence, the required assertion holds.

 Step 4\quad Prove that
$\sup_{h\in\mathscr{H}}\inf_{x\in(0,
D)}R(h)(x)\geqslant\lambda_{p}$ when $u$ and $v'$ are continuous.

Noticing that $$\hat\nu(x, D)\bigg(\int_0^xf^{p-1}\d\mu\bigg)^{p^*-1}\leqslant fI\!I(f)(x)^{p^*-1}\leqslant\hat\nu(x, D)\bigg(\int_0^Df^{p-1}\d\mu\bigg)^{p^*-1},$$
If $\hat\nu(0, D)<\infty$, then choose  $f\in L^{p-1}(\mu)$ to be a positive function  such that $g=fI\!I(f)^{p^*-1}<\infty$.
 Set $\bar h=g'/g$. Then $\bar h\in\mathscr{H}$ since $u$ and $v'$ are continuous.
 Moreover, $L_pg=-u f^{p-1}$ and
$$R(\bar h)=-\frac{1}{u g^{p-1}}L_pg=\frac{f^{p-1}}{g^{p-1}}>0.$$
If $\hat\nu(0, D)=\infty$, then set $\bar h=0$. So $R(\bar h)=0$.
In other words, we always have  $$\sup_{h\in\mathscr{H}}\inf_{x\in(0,D)}R(h)(x)\geqslant0.$$

Without loss of generality, assume that $\lambda_p>0$ and $g$ is an eigenfunction of $\lambda_p$, i.e.,
$$L_pg=-\lambda_pu  |g|^{p-2}g.$$
Let $h=g'/g\in \mathscr{H}$. Then $R(h)=\lambda_p$ and the assertion holds.

Step 5\quad Prove that the supremum in the lower estimates can be attained.

Since
$$0=\lambda_p\geqslant\inf_{x\in (0, D)}I\!I(f)(x)^{-1}\geqslant0,\qquad
0=\lambda_p\geqslant\inf_{x\in (0, D)}I(f)(x)^{-1}\geqslant0$$
for every $f$ in the set defining $\lambda_p$, the assertion is
clear for the case that $\lambda_p=0$. Similarly, the conclusion
holds for operator $R$ as seen from the preceding proof in Step 4.
For the case that $\lambda_p>0$, assume that $g$  is an
eigenfunction corresponding to $\lambda_p$. Let $\bar h=g'/g\in
\scr{H}$. Then $R(\bar h)=\lambda_p$, $I(g)^{-1}\equiv\lambda_p$ by
letting $f=\lambda_p^{p^*-1}g$ in \eqref{C-N-D} and
$I\!I(g)^{-1}\equiv\lambda_p$ by letting $f=\lambda_p^{p^*-1}g$ in
\eqref{C-Nf1} whenever $g(D)=0$.

Now, it remains to show that the vanishing property of eigenfunction
at $D$, which is proved in the following proposition by  using the
variational formula  proved in Step 1 above.

\begin{prop}\label{CNp1} Let $g$ be an  a.e. eigenfunction of $\lambda_p>0$.
Then $g(D)=0$.
\end{prop}
{\it Proof }\quad Let $f=g-g(D)$. Then $f\in\mathscr{F}_{I\!I}$. By
\eqref{C-Nf1}, we have
$$f(x)=\lambda_p^{p^*-1}\int_x^D\hat{v}(t)\bigg(\int_0^tg^{p-1}\d\mu\bigg)^{p^*-1}\d
t.$$
We prove the proposition by dividing it into two cases. Denoted
by
$$M(x)=\int_x^D\hat{v}(t)\bigg(\int_0^t\d\mu\bigg)^{p^*-1}\d t.$$

(a) If $M(x)=\infty$, then $f(x)=g(x)-g(D)<\infty$ and
$$\aligned\lambda_p^{1-p^*}f(x)&=\int_x^D\hat{v}(t)\bigg(\int_0^tg^{p-1}\d\mu\bigg)^{p^*-1}\d t
>g(D)M(x)=\infty\endaligned$$
once $g(D)\neq0$. So there is a contradiction.

(b) If $M(x)<\infty$, then
$$fI\!I(f)(x)^{p^*-1}=\int_x^D\hat{v}(t)\bigg(\int_0^t(g-g(D))^{p-1}\d\mu\bigg)^{p^*-1}\d t
<g(0)M(0)<\infty.$$
Replacing $f$ in the denominator of $I\!I(f)$ with this term and
using  Cauchy's mean-value theorem  twice, we have
 $$\aligned
\sup_{(0, D)}I\!I(f)&\leqslant\!
\frac{1}{\lambda_p}\sup_{(0,D)}\frac{f^{p-1}}{g^{p-1}}=
\frac{1}{\lambda_p}\!\sup_{x\in(0,D)}\!\bigg(1\!-\!\frac{g(D)}{g(x)}\bigg)^{p-1}\!\!
=\frac{1}{\lambda_p}\!\bigg(1\!-\!\frac{g(D)}{g(0)}\bigg)^{p-1}.\endaligned$$
The last equality comes from the fact that $g\downdownarrows$.
If $g(D)>0$, then
$$\lambda_p^{-1}\leqslant\inf_{f\in\scr{F}_{I\!I}}\sup_{x\in(0,
D)}I\!I(f)(x)\leqslant\sup_{x\in(0, D)}I\!I(f)(x)<\lambda_p^{-1},$$
which is a contradiction. Therefore, we must have
$g(D)=0$.$\qquad\Box$

By now, we have finished the proof of the lower estimates of
$\lambda_p$. Dually, one can prove the upper estimates without too
much difficulty.
We ignore the details here.

The following lemma or its variants have been used many times before
(cf., \rf{Chen2010}{Proof of Theorem 3.1}, \rf{Chen2005}{page 97},
or \cite{JM2012}, and the earlier publications therein). It is
essentially an application of the integration by parts formula, and is
a key to the proof of Theorem \ref{CNth2}.
\begin{lem}\label{CNL1}
 Assume that $m$ and $ n$ are two non-negative locally integrable functions. For $p>1$, define
 $$S(x)=\bigg(\int_{x}^{D}n(y)\d y\bigg)^{p-1},\qquad M(x)=\int_{0}^{x}m(y)\d y$$
 and $c_{0}=\sup_{x\in(0,D)}S(x)M(x)<\infty$.
 Then
  $$\int_{0}^{x}m(y)S(y)^{{p^*r}/{p}}\d y\leqslant \frac{c_{0}}{1-{p^*r}/{p}}S(x)^{({p^*r}/{p})-1},
  \qquad r\in(0, {p}/{p^*}).
 $$
 \end{lem}

 {\noindent\it Proof of Theorem $\ref{CNth2}$}
\quad First, we prove that $\lambda_p\geqslant\big(k(p)\sigma_p\big)^{-1}$. Fixing $r\in(0, p/p^*)$, let
$f(x)=\hat{\nu}(x, D)^{p^*r/p}$.
 Applying  $m(x)=u(x),\ n(x)=\hat{v}(x)$ to Lemma \ref{CNL1}, we have $M(x)=\mu(0,x)$,
 $S(x)=\hat{\nu}(x, D)^{p-1}$, $c_{0}=\sigma_p$ and
$$
\int_{0}^{x}\hat{\nu}(y, D)^{r}\mu(\d y)\leqslant\frac{\sigma_p}{1-{p^*r}/{p}}
\hat{\nu}(x, D)^{r-(p/p^*)}.$$
Since
 $$|f'|^{p-2}f'=-\bigg(\frac{p^*r}{p}\hat{\nu}(\cdot,
 D)^{(p^*r/p)-1}\hat{v}(\cdot)\bigg)^{p-1},
$$
we have
\begin{equation}\label{NDstar}\sup_{x\in (0, D)} I(f)(x)\leqslant
\frac{[{p}/(p^*r)]^{p-1}}{1-{p^*r}/{p}}\sigma_p.
\end{equation}
By Theorem \ref{CNth1}\,(1),  \eqref{NDstar}, and an optimization
with respect to $r\in(0, p/p^*)$, we obtain
$$
\lambda_{p}^{-1}\leqslant\Big(\sup_{f\in\mathscr{F}_I}\inf_{x\in(0,
D)} I(f)(x)^{-1}\Big)^{-1}\leqslant
p{p^*}^{p-1}\sigma_p=k(p)\sigma_p.$$

Now we prove that $\lambda_p\leqslant\sigma_p^{-1}$. For fixed
$x_{0},\, x_{1}\in(0,D)$ with $x_{0}<x_{1}$, let $f(x)=\hat{\nu}(x\vee x_{0}, D)
\mathbbold{1}_{[0,x_{1})}(x)$.
Then
$$I(f)(x)=\hat{\nu}(x_0, D)^{p-1}\mu(0,x_{0})+\int_{x_{0}}^{x}\hat{\nu}(t, D)^{p-1}\mu(\d t),
\qquad x\in(x_{0},x_{1})$$
and $I(f)(x)=\infty$ on $[0, x_0]\cup[x_1, D]$ by convention $
{1}/{0}=\infty$. Combining with Theorem \ref{CNth1}\,(1), we have
 $$\lambda_{p}^{-1}\geqslant\inf_{x< x_{1}}I(f)(x)=\hat{\nu}(x_0, D)^{p-1}\mu(0,x_{0}),
 \qquad x_0<x_1.$$
Thereby the assertion that $\lambda_p\leqslant\sigma_p^{-1}$ follows by letting
 $x_{1}\rightarrow D$. Since $$\mu(0,x)^{p^*-1}\hat\nu(x, D)\leqslant\int_x^D\mu(0,s)^{p^*-1}\hat v(s)\d s\leqslant \int_0^D\mu(0,s)^{p^*-1}\hat v(s)\d s,$$
 the assertions hold.
 $\qquad\Box$

 From the proof above, it is easy to understand why we choose the test function as
 $f=\hat{\nu}[\cdot, D]^{1/p^*}$ in \rf{CWZ2013-1}{Proof of Theorem 2.3 (a)} in the discrete case.

\medskip
\noindent{\it Proof of Theorem $\ref{CNth3}$}
 \quad Using  Cauchy's mean-value theorem  and definitions of $\delta_n'$, $\delta_n$,
 $\bar\delta_n$ and $\lambda_p$,
 it is not hard to show the most of the results
except that $\bar\delta_{n+1}\geqslant\delta_n'$. Put $f=f_n^{x_0, x_1}$ and $g=f_{n+1}^{x_0, x_1}$.
Then $g=fI\!I(f)^{p^*-1}$. By simple calculation, we have
$$\begin{aligned}
   D_p(g)&=\int_{0}^{x_{1}}|g'|^{p-1}|g'|(x)v(x)\d x
=\int_{0}^{x_{1}}v(x)^{-1}\int_{0}^{x}f^{p-1}\d\mu|g'(x)|v(x)\d x
\end{aligned}$$ Exchanging the order of the integrals, we have
$$\begin{aligned}D_p(g)&=-\int_{0}^{x_1}f^{p-1}(t)\mu(\d t)\int_{t}^{x_1}g'(x)\d x
\quad(\text{by Fubini's Theorem})\\
   &\leqslant\int_{0}^{x_1}f^{p-1}(t)g(t)\mu(\d t)\quad(\text{since } g(x_1)\geqslant0)\\
&\leqslant\int_{0}^{x_1}g^p\d\mu\sup_{t\in(0, x_1)}\bigg(\frac{f(t)}{g(t)}\bigg)^{p-1}\\
&\leqslant \mu(|g|^p)\sup_{x\in(0, x_1)}I\!I(f)(x)^{-1}.
 \end{aligned}$$
So the required assertion holds. $\qquad\Box$
\medskip

\noindent{\it Proof of Corollary $\ref{CNcor1}$ }\quad
 (a) \; The calculation of $\delta_1$ is simple. We compute $\delta_1'$ first.  Consider the term
 $\inf_{x< x_1}I\!I(f_1^{x_0,x_1})(x)$.
By  calculation, we obtain that for $x\in(x_0, x_1)$, the numerator
of $\big(I\!I(f_1^{x_0,x_1})(x)^{p^*-1}\big)'\big|_x$  equals
$$\hat{v}(x)\bigg[\int_{x}^{x_1}\hat{v}(s)\bigg(\int_0^s(f_1^{x_0,x_1})^{p-1}\d\mu\bigg)^{p^*-1}\d s
-\hat{\nu}(x,x_1)\bigg(\int_0^x(f_1^{x_0,x_1})^{p-1}\d\mu\bigg)^{p^*-1}\bigg]$$
which is obviously non-negative. So
$$I\!I(f_1^{x_0,x_1})(x)=\bigg[\frac{1}{\hat{\nu}(x,x_1)}\int_x^{x_1}\hat{v}(s)
\bigg(\int_0^s\hat{\nu}(t\vee x_0, x_1)^{p-1}\mu(\d t)\bigg)^{p^*-1}\d s\bigg]^{p-1}$$ is
increasing in $x\in(x_0, x_1)$.
Hence,
 $$\aligned\delta_1'&=\sup_{x_0<x_1}\bigg[\frac{1}{\hat{\nu}(x_0,x_1)}\int_{x_0}^{x_1}\hat{v}(s)
 \bigg(\int_0^s\hat{\nu}(t\vee x_0, x_1)^{p-1}\mu(\d t)\bigg)^{p^*-1}\d s\bigg]^{p-1}\\
&=\sup_{x_0\in (0, D)}\frac{1}{\hat{\nu}(x_0,
D)^{p-1}}\bigg[\int_{x_0}^{D}\hat{v}(s)\bigg(\int_0^s\hat{\nu}(t\vee x_0,
D)^{p-1}\mu(\d t)\bigg)^{p^*-1}\d s\bigg]^{p-1}.
\endaligned$$
In the last equality, we have used the fact that
$I\!I(f_1^{x_0,x_1})(x_0)$ is increasing in $x_1\in [x_0, D]$. Indeed, let
$$N_k(s,y)=\int_{x_0}^{s}\hat{\nu}(t, y)^{k}\mu(\d t),\quad f(s, y)=\hat{v}(s)N_{p-1}(s,y)^{p^*-1}.$$
Then
 $$\aligned
 I\!I(f_1^{(x_0, y)})(x_0)^{p^*-1}&\!\!=\!\frac{1}{\hat{\nu}(x_0, y)}\bigg[\int_{x_0}^yf(s, y)\d s
 \!+\!\!\int_{x_0}^y\hat{v}(s)\d s\int_0^{x_0}\hat{\nu}(x_0, y)^{p-1}\mu(\d t)\bigg]\\
 &=\frac{1}{\hat{\nu}(x_0, y)}\int_{x_0}^yf(s, y)\d s+\mu(0, x_0)\hat{\nu}(x_0, y)^{p-1}\\
 &=:H_1(y)+H_2(y),
 \endaligned$$
 and
 $$\aligned \frac{\partial}{\partial y} N_{p-1}(s,y)&=\int_{x_0}^s(p-1)\hat{\nu}(t,y)^{p-2}\hat{v}(y)
 \mu(\d t)\\
 &=(p-1)\hat{v}(y)N_{p-2}(s,y);\endaligned$$
 $$\aligned
  \frac{\partial}{\partial y} f(s, y)&=(p^*-1)\hat{v}(s)N_{p-1}(s,y)^{p^*-2}\frac{\partial }
  {\partial y}N_{p-1}(s,y);\endaligned$$
 $$\aligned
\frac{\partial }{\partial y}\int_{x_0}^yf(s, y)\d
s&=\int_{x_0}^y\frac{\partial}{\partial y} f(s, y)\d
s+f(y,y).\endaligned$$ Hence, the numerator of $\d H_1/\d y $ equals
$$\aligned&\bigg(\frac{\partial }{\partial y}\int_{x_0}^yf(s, y)\d s\bigg)
\hat{\nu}(x_0, y)-\hat{v}(y)\int_{x_0}^yf(s,y)\d s\\
&=\hat{\nu}(x_0, y)\hat{v}(y)\int_{x_0}^y\hat{v}(s)N_{p-1}(s,y)^{p^*-2}N_{p-2}(s,y)\d s\\
&\hskip3cm+\hat{\nu}(x_0, y)f(y,y)-\hat{v}(y)\int_{x_0}^yf(s, y)\d s\\
&=\hat{v}(y)\bigg(\hat{\nu}(x_0, y)\int_{x_0}^y\hat{v}(s)N_{p-1}(s,y)^{p^*-2}N_{p-2}(s,y)\d s\\
&\hskip3cm-
\int_{x_0}^y\hat{v}(s)N_{p-1}(s,y)^{p^*-1}\d s
\bigg)+\hat{\nu}(x_0, y)f(y,y).\endaligned$$
Since $\hat{\nu}(x_0, y)N_{p-2}(s, y)-N_{p-1}(s,y)\geqslant0$ for $s\in[x_0, y]$, we see that
$\d H_1/\d y$  is positive. It is obvious that $\d H_2/\d y$ is positive. So
$I\!I(f_1^{x_0, y})(x_0)$ is increasing in $y$ and the required assertion holds.

(b)\; Compute $\bar{\delta}_1$. By  definition of $\bar{\delta}_1$,
we have
 $$\begin{aligned}
 \|f_1^{x_{0}, x_{1}}\|_p^{p}&=\int_{0}^{x_{1}}\bigg(\int_{x_{0}\vee x}^{x_{1}}\hat{v}(s)\d s\bigg)^{p}
\mu(\d x)\\
&
=\mu(0,x_{0})\hat{\nu}(x_{0},x_{1})^{p}+\int_{x_{0}}^{x_{1}}
\bigg(\int_{x}^{x_{1}}\hat{v}(t)\d t\bigg)^{p}\mu(\d x),
\\D_p(f_1^{x_{0}, x_{1}})&=\int_{x_{0}}^{x_{1}}\hat{v}(t)^pv(t)\d t=\hat{\nu}(x_0,x_{1}).\end{aligned}$$
 Hence,
$$\begin{aligned}
\bar{\delta}_1&=\sup_{x_{0}<x_{1}}\bigg(\mu(0,x_{0})\hat{\nu}(x_0,x_{1})^{p-1}+
\frac{1}{\hat{\nu}(x_0,x_{1})}\int_{x_{0}}^{x_{1}}\hat{\nu}(s, x_1)^{p}\mu(\d s)\bigg)
\\&=\sup_{x_{0}\in(0, D)}\bigg(\mu(0,x_{0})\hat{\nu}(x_0,D)^{p-1}+
\frac{1}{\hat{\nu}(x_0, D)}\int_{x_{0}}^D\hat{\nu}(s, D)^{p}\mu(\d s)\bigg)
\end{aligned}$$
In the second equality, we have used the fact that:
$$\mu(0,x_{0})\hat{\nu}(x_0,x_{1})^{p-1}+
\frac{1}{\hat{\nu}(x_0,x_{1})}\int_{x_{0}}^{x_{1}}\hat{\nu}(s, x_1)^p\mu(\d s)
\quad\uparrow\quad\text{ in } x_1.$$
Indeed, it suffices to show that
$$\frac{1}{\hat{\nu}(x_0,x)}\int_{x_{0}}^{x}\hat{\nu}(s, x)^p\mu(\d s)
\leqslant\frac{1}{\hat{\nu}(x_0, y)}\int_{x_{0}}^y \hat{\nu}(s, y)^p\mu(\d s),\qquad x_0\leqslant x<y,$$
which is equivalent to
$$\frac{1}{\hat{\nu}(x_0,y)}\int_x^y \hat{\nu}(s, y)^p\mu(\d s)
+\int_{x_{0}}^{x}\frac{\hat{\nu}(s, y)^p}{\hat{\nu}(x_0,y)}
-\frac{\hat{\nu}(s, x)^p}{\hat{\nu}(x_0,x)}\mu(\d s)\geqslant0.$$
Since $p>1$ and $\hat{\nu}(t, x)\leqslant \hat{\nu}(x_0, x)$ for $x\geqslant t\geqslant x_0$, we have
$$\frac{\hat{\nu}(t, y)^p}{\hat{\nu}(t,x)^p}=\bigg[\frac{\hat{\nu}(t, x)
+\hat{\nu}(x, y)}{\hat{\nu}(t, x)}\bigg]^p\geqslant1+\frac{\hat{\nu}(x, y)}{\hat{\nu}(t, x)}\geqslant1
+\frac{\hat{\nu}(x, y)}{\hat{\nu}(x_0, x)}=\frac{\hat{\nu}(x_0, y)}{\hat{\nu}(x_0, x)}$$
for $ t\geqslant x_0$
and the required assertion holds.

(c)\; Comparing $\delta_1'$ and $\bar\delta_1$. It is easy to see that
$$\aligned
\int_x^D\hat{\nu}(\cdot, D)^p\d\mu&=\int_x^D\hat{\nu}(t, D)^{p-1}\int_t^D\hat{v}(s)\d s\mu(\d t)\\
&=\int_x^D\hat{v}(s)\int_x^s\hat{\nu}(t, D)^{p-1}\mu(\d t)\d s ;
\endaligned$$
$$\aligned
&\mu(0, x)\hat{\nu}(x, D)^{p}=\int_x^D\hat{v}(s)\int_0^x\hat{\nu}(x, D)^{p-1}\mu(\d t)\d s.
\endaligned$$
Let $a_x(s)=\hat{v}(s)\big/\hat{\nu}(x, D)$ for $s\in (x, D)$. Noticing that $a_x$ is a probability
on $(x, D)$, by the increasing property of moments $\mathbb{E}(|X|^s)^{1/s}$ in $s>0$ and combining
the preceding assertions (a) and (b),  we have
$$\aligned
\bar\delta_1&=\sup_{x\in(0, D)}
\int_x^Da_x(s)\int_0^s\hat{\nu}(t\vee x, D)^{p-1}\mu(\d t)\d s\\
&\leqslant\sup_{x\in(0, D)}
\bigg[\int_x^Da_x(s)\bigg(\int_0^s\hat{\nu}(t\vee x, D)^{p-1}\mu(\d t)\bigg)^{p^*-1}\d s\bigg]^{p-1}\;
(\text{if } p^*-1>1)\\
&=\delta_1'.
\endaligned$$
Similarly, if $p^*-1<1$ (i.e., $p>2$), then $\bar\delta_1\geqslant\delta_1'$.

(d)\; Prove that $\bar\delta_1\leqslant p\sigma_p$. Using the
integration by parts formula, we have
$$\begin{aligned}
 \int_{x_0}^x\hat{\nu}(y, D)^p\mu(\d y)
 &=\hat{\nu}(y, D)^p\mu(0,y)\big|_{x_0}^x+p\int_{x_0}^{x}\hat{\nu}(y, D)^{p-1}\hat{v}(y) \mu(0,y)\d y\\
&\leqslant \sigma_p{\hat{\nu}(x, D)}-\hat{\nu}(x_0, D)^p\mu(0,x_0)+p\sigma_p\int_{x_0}^x\hat{v}(y)\d y
\end{aligned}$$
Since $\hat{\nu}(x, D)<\infty$, letting $x\to D$, we have
$$\aligned\bar\delta_1&=\sup_{x_{0}\in(0, D)}\bigg(\mu(0,x_{0})\hat{\nu}(x_0, D)^{p-1}+
\frac{1}{\hat{\nu}(x_0, D)}\int_{x_{0}}^D\hat{\nu}(\cdot, D)^p\d\mu\bigg)\\
&\leqslant\sup_{x_{0}\in(0, D)}\bigg[\mu(0,x_{0})\hat{\nu}(x_0, D)^{p-1}+\\
&\hskip2cm
\frac{1}{\hat{\nu}(x_0, D)}\bigg(-\hat{\nu}(x_0, D)^p\mu(0,x_0)+p\sigma_p\int_{x_0}^D\hat{v}(y)\d y
\bigg)\bigg]\\
&=p\sigma_p,\endaligned$$ and the required assertion
holds.$\qquad\Box$
\noindent\\[4mm]

\setcounter{section}{4}
\setcounter{thm}{0}

\noindent{\bbb 4\quad DN-case}\\[0.1cm]
\noindent From  now on, we concern on $p\,$-Laplacian eigenvalue with
DN-boundaries. We use the same notation as the previous ND-case
since they play the similar role but have different meaning in
different context. Let $D\leqslant\infty$, $p>1$. The p-Laplacian
eigenvalue problem with DN-boundary conditions is
\begin{equation}\begin{cases}
\text{Eigenequation}: &L_pg(x)=-\lambda u(x)|g|^{p-2}g(x);\\
\text{DN-boundaries}: &g(0)=0,\quad g'(D)=0\text{ if } D<\infty \label{CDN-eigen}\end{cases}
\end{equation}
The first eigenvalue $\lambda_p$ has the following classical variational formula:
 \begin{equation}\label{DN-V}\lambda_p\!=\!\inf\!\bigg\{\frac{D_p(f)}{\mu(|f|^p)}:
 f(0)\!\!=\!0,\; f\!\!\ne\! 0, f\!\in\!\scr{C}[0, D], v^{p^*-1}f'\!\in\!\scr{C}(0, D), D_p(f)<\infty\bigg\}.\end{equation}
Correspondingly, we are also estimating the optimal constant
$A:=\lambda_p^{-1}$ in the {\it weighted Hardy inequality}:
\begin{equation*}\mu(|f|^{p})\leqslant A D_{p}(f),\quad f(0)=0,\quad f\in\mathscr{D}(D_{p}).\end{equation*}
For $p>1$, define $\hat{v}=v^{1-p^*}$ and $\hat{\nu}(\d
x)=\hat{v}(x)\d x$. We use the following operators:
$$\aligned
I(f)(x)&=\frac{1}{\big(v f'|f'|^{p-2}\big)(x)}\int_x^Df^{p-1}\d\mu
\quad(\text{single integral form})\\ I\!I(f)(x)&=\frac{1}{f^{p-1}(x)}\bigg[\int_0^x\hat{v}(s)
\bigg(\int_s^Df^{p-1}\d\mu\bigg)^{p^*-1}\d s\bigg]^{p-1}\\
&\hskip7cm\quad(\text{double integral form})\\R(h)(x)&=-u^{-1}\big\{|h|^{p-2}\big[v'h+(p-1)v(h^2+h')
\big]\big\}(x) \quad(\text{differential form}).
\endaligned$$
The three operators above have domains, respectively, as follows.
$$\aligned&\mathscr{F}_{I}=\{f\in\mathscr{C}[0,D]: v^{p*-1}f'\in\mathscr{C}(0,D), f(0)=0\text{ and } f'|_{(0,D)}>0 \},\\
&\mathscr{F}_{I\!I}=\{f:f\in\mathscr{C}[0,D], f(0)=0\text{ and } f|_{(0, D)}>0\}.\\
&\mathscr{H}=\bigg\{h: h\in\mathscr{C}^{1}(0,
D)\cap\mathscr{C}[0,D], h|_{(0,D)}>0 \text{ and } \int_{0+}h(u)\d u=\infty\bigg\},
\endaligned$$
where $\int_{0+}$ means $\int_0^{\varepsilon}$ for sufficiently small $\varepsilon>0$.
Some modifications are needed when studying the upper estimates.
$$\aligned&\widetilde{\mathscr{F}}_{I}=\big\{f\in\mathscr{C}[0, x_0]: f(0)=\!0, v^{p*-1}f'\in
\mathscr{C}(0, x_0), f'|_{(0, x_0)}>0 \text{ for some } \\
&\hskip1.7cm x_{0}\in(0, D),
\text{ and }f=f(\cdot\wedge x_{0})\big\},\\
&\widetilde{\mathscr{F}}_{I\!I}=\big\{f: f(0)\!=\!0, \exists\, x_0\!\in\!(0, D)
\text{ such that }f\!\!=\!\!f(\cdot\wedge x_0)\!>\!0\text{ and } f\!\in\!\mathscr{ C}[0, x_{0}]\big\},\\
&\widetilde{\mathscr{H}}=\bigg\{h: \exists\, x_{0}\in(0,
D)\text{ such that } h\in\mathscr{C}[0,x_{0}]\cap\mathscr{C}^{1}(0, x_{0}), h|_{(0, x_{0})}>0, \\
&\hskip 1.5cm h|_{[x_0, D]}\!=\!0, \int_{0+}\!h(u)\d
u\!=\!\infty,\text{ and\,} \sup_{(0,
x_0)}\!\big[v'h+(p-1)(h^2+h')v\big]\!<\!0\bigg\}.\endaligned$$ When
$D=\infty$,  replace $[0, D]$ and $(0, D]$ with $[0, D)$ and $(0,
D)$, respectively. Besides, we also need the following notation:
$$\widetilde{\mathscr{F}}_{I\!I}'=\big\{f: f(0)=0,  f\in\mathscr{ C}[0, D]
\text{ and }fI\!I(f)\in L^p(\mu) \big\}.$$

If $\mu(0,D)=\infty$, then $\lambda_p$ defined by \eqref{DN-V} is trivial.
Indeed, let $$f=\mathbbold{1}_{(\delta,\, D]}+h\mathbbold{1}_{[0, \delta]},$$
 where $h$ is chosen such that $h(0)=0$ and $f\in {\scr C}^1(0,D)\cap {\scr C}[0,D]$ (for example, $h(x)=-x^{2}\cdot\delta^{-2}+2x\cdot\delta^{-1}$).
Then $D_p(f)\in (0,\infty)$ and $\mu\big(|f|^p\big)=\infty$.
It follows that  $\lambda_p=0.$

Otherwise, $\mu(0,D)<\infty$. Then for every $f$ with $\mu\big(|f|^p\big)=\infty,$ by setting $f^{(x_{0})}=f(\cdot\wedge x_{0})\in L^p(\mu)$,  we have
$$
\aligned\infty>D\big(f^{(x_{0})}\big)\to D(f),\quad  \infty>\mu\Big({|f^{(x_{0})}|}^p\Big)\rightarrow
\mu\big(|f|^p\big)\quad\text{as}\; x_{0}\rightarrow D.
\endaligned$$
 In other words, for  $f\notin L^p(\mu)$,  both $\mu\big(|f|^p\big)$ and $D_p(f)$ can be approximated by a sequence of functions belonging to $L^p(\mu)$. Hence,
we can rewrite $\lambda_p$ as follows.
\begin{equation}\label{DN-V1}
\lambda_p=\inf\big\{D_p(f):  \mu\big(|f|^p\big)=1, f(0)=0, \text{ and } f\in {\scr C}^1(0,D)\cap {\scr C}[0,D]\big\}.
\end{equation}
In this case,  we also have
 $$\aligned
\lambda_p=\inf\big\{ &D_p(f):\mu\big(|f|^p\big)=1,\, f(0)=0,\,f=f(\cdot\wedge x_{0}), \\&f\in {\scr C}^1(0,x_{0})\cap {\scr C}[0,x_{0}]\text{ for } \text{some } x_{0}\in (0,D)\big\}.\endaligned$$
We are now ready to state the main results in the present context.
\begin{thm}
\label{CDth1} Assume that $\mu(0, D)<\infty$. For $p>1$, the following variational formulas hold for $\lambda_p$ defined by \eqref{DN-V1} $($equivalently, \eqref{DN-V}$)$.
\begin{itemize}\setlength{\itemsep}{-0.8ex}
\item[$(1)$] Single integral forms:
$$
\aligned
\inf_{f\in\widetilde{\mathscr{F}_I}}\sup_{x\in(0,D)}I(f)(x)^{-1}
=\lambda_p=\sup_{f\in\mathscr F_{I}}\,\inf_{x\in(0,D)}I(f)(x)^{-1},
\endaligned$$
\item[$(2)$]Double integral forms:
$$\aligned
&\lambda_p=\inf_{f\in \widetilde{\mathscr{F}}_{I}}\sup_{x\in(0,
D)}I\!I(f)(x)^{-1} =\inf_{f\in
\widetilde{\mathscr{F}}_{I\!I}\cup\widetilde{\mathscr{F}}_{I\!I}'}\sup_{x\in(0,
D)}I\!I(f)(x)^{-1},\\
&\lambda_p=\sup_{f\in \mathscr{F}_{I}} \inf_{x\in(0, D)}I\!I(f)(x)^{-1}=
\sup_{f\in \mathscr{F}_{I\!I}} \inf_{x\in(0, D)}I\!I(f)(x)^{-1}
.
\endaligned$$
\end{itemize}
Moreover, if $u$ and $v'$ are continuous,  then we have additionally
\begin{itemize}\setlength{\itemsep}{-0.8ex}
\item[$(3)$]  differential forms:
$$
\aligned
 \inf_{h\in\widetilde{\mathscr{H}}}\,\sup_{x\in(0,
D)}R(h)(x)=\lambda_p=\sup_{h\in\mathscr{H}}\,\inf_{x\in(0,
D)}R(h)(x).
\endaligned$$
\end{itemize}\end{thm}

Define $k(p)=p{p^*}^{p-1}$ and
$$\sigma_p=\sup_{x\in(0,D)}\mu(x, D)\hat{\nu}(0, x)^{p-1}.$$
 As an application of the variational formulas in Theorem \ref{CDth1} (1), we have the following theorem
 which was also known in 1990's (cf. \rf{BA1990}{Lemmas 3.2 and 3.4 on pages 22 and 25, respectively}).
\begin{thm}
\label{CDth2}$(\text{\rm Criterion and basic estimates})$ For $p>1$, $\lambda_p>0$ if and only if
$\sigma_p<\infty$.
Moreover,
 \begin{equation*}
 (k(p)\sigma_p)^{-1}\leqslant\lambda_{p}\leqslant\sigma_p^{-1}.
 \end{equation*}
 In particular, we have $\lambda_p=0$ if $\mu(0, D)=\infty$ and $\lambda_p>0$ if
 $$\int_0^D\mu(s, D)^{p^*-1}{\hat\nu}(\d s)<\infty.$$
 \end{thm}

The next result is an application of the variational formulas in
Theorem \ref{CDth1} (2).

 \begin{thm}\label{CDth3}$(\text{\rm Approximating procedure})$
Assume that $\mu(0, D)<\infty$ and $\sigma_p\!<\!\infty$.
\begin{itemize}
\item[$(1)$] Let $f_1=\hat{\nu}(0, \cdot)^{1/p^*}$,
$f_{n+1}=f_{n}I\!I(f_{n})^{p^*-1}$ and $\delta_{n}=\sup_{x\in(0,D)}I\!I(f_{n})(x)$ for  $n\geqslant1$.
 Then $\delta_{n}$ is decreasing in $n$ and
$$\lambda_{p}\geqslant\delta_{n}^{-1}\geqslant(k(p)\sigma_p)^{-1}.$$
\item[$(2)$] For fixed $x_{0}\in(0,D)$,  let
$$f_1^{(x_{0})}=\hat{\nu}(0, \cdot\wedge x_0),\qquad
f_{n}^{(x_{0})}=f_{n-1}^{(x_0)}
I\!I\big(f_{n-1}^{(x_0)}\big)(\cdot\wedge x_0)^{p^*-1}$$ and
$\delta_{n}'=\sup_{x_{0}\in(0, D)}\inf_{x\in(0,D)}
I\!I\big(f_{n}^{(x_{0})}\big)(x)$ for $n\geqslant1$. Then
$\delta_{n}'$ is increasing in $n$ and
$$\sigma_p^{-1}\geqslant{\delta_{n}'}^{-1}\geqslant\lambda_{p}.$$
Moreover, define
$$\bar{\delta}_n=\sup_{x_{0}\in(0,D)}\frac{\big\|f_n^{(x_{0})}\big\|_p^p}{D_{p}\big(f_{n}^{(x_{0})}\big)},
\qquad n\geqslant 1.$$
Then $\bar{\delta}_n^{-1}\geqslant\lambda_{p}$ and
$\bar{\delta}_{n+1}\geqslant\delta_n'$ for $n\geqslant1$.
\end{itemize}
\end{thm}
Most of the result in Corollary \ref{CDcor1} below can be obtained directly from Theorem
\ref{CDth3}.
\begin{cor}
\label{CDcor1}$(\text{\rm Improved estimates})$  Assume that $\mu(0, D)<\infty$ and $\lambda_p>0$.
We have $$\sigma_p^{-1}\geqslant\delta_1'^{-1}\geqslant\lambda_p\geqslant\delta_1^{-1}
\geqslant(k(p)\sigma_p)^{-1},$$
where
$$\aligned\delta_1&=\sup_{x\in(0, D)}\bigg[\frac{1}{\hat{\nu}(0, x)^{1/p^*}}\int_0^x\hat{v}(s)
\bigg(\int_s^D\hat{\nu}(0, t)^{p/{p^*}^2}\mu(\d t)\bigg)^{p^*-1}\d s\bigg]^{p-1};\endaligned$$
$$\aligned\delta_1'&=\sup_{x\in (0, D)}\frac{1}{\hat{\nu}(0, x)^{p-1}}\bigg[\int_0^{x}\hat{v}(s)
\bigg(\int_s^D\hat{\nu}(0, t\wedge x)^{p-1}\mu(\d t)\bigg)^{p^*-1}\d s\bigg]^{p-1}.\endaligned$$
Moreover, $$\aligned
\bar\delta_1&\!=\!\sup_{x\in(0, D)}\bigg(\mu(x, D)\hat{\nu}(0, x)^{p-1}+
\frac{1}{\hat{\nu}(0, x)}\int_0^{x}\hat{\nu}(0, t)^p\mu(\d t)\bigg)\in[\sigma_p, p \sigma_p],
\endaligned$$
and
$\bar\delta_1\geqslant\delta_1$ for $p\geqslant2$ and $\bar\delta_1\leqslant\delta_1'$
for $1<p\leqslant2$.
\end{cor}
When $p=2$, the equality $\delta_1=\bar\delta_1$ was proved in
\rf{CWZ2013}{Theorem 6}.

Most of the results in this section are parallel to that in Section
2. One may follow Section 3 or \cite{CWZ2013,JM2012} to complete the
proofs without too many difficulties. The details are omitted here.
Instead, we prove some properties of the eigenfunction $g$, which
are used in choosing the test functions for the operators. 
\begin{lem}Let $(\lambda_p, g)$ be a solution to \eqref{CDN-eigen}, $g\ne 0$. Then $g'$ does
not change sign, and so does $g$.
\end{lem}
{\noindent\it Proof}\quad First,  the solution provided by Lemma \ref{Cexists} is trivial:
$g=0$, if the given constants
$A$ and $B$ are zero. Because we are
in the situation that $g(0)=0$, we can assume that $g'(0)\ne 0$. Next, we prove that $g'$ dose not change
sign by seeking a contradiction. If there exists $x_0\in(0, D)$ such
that $g'(x_0)=0$, then $g(x_0)\neq0$ by \rf{P2006}{Lemma 2.3}. Let $\bar g=g\mathbbold{1}_{[0, x_0]}+g(x_0)\mathbbold{1}_{(x_0, D]}$. By simple calculation, we obtain
$$D_p(\bar g)=(-L_p\bar g, \bar g)_{\mu}=\lambda_p\mu_{0, x_0}(|g|^p).$$
So
$$\lambda_p\leqslant\frac{D_p(\bar g)}{\mu(|\bar g|^p)}=\frac{\lambda_p\mu_{0, x_0}(|g|^p)}{\mu_{0,x_0}
\big(|g|^p\big)+\mu(x_0, D)|g(x_0)|^p}<\lambda_p,$$
which is a contradiction. Therefore $g'$ does not change sign. Since
$g(0)=0$, the second assertion holds naturally.
\noindent\\[4mm]

\noindent\bf{\footnotesize Acknowledgements}\quad\rm
{\footnotesize  The work is supported in part by NSFC (Grant No.11131003), SRFDP
 (Grant No. 20100003110005), the ``985'' project from the Ministry of Education in China and the
 Fundamental Research Funds for the Central Universities.
The authors also thank Professor Yong-Hua Mao for his helpful
comments and suggestions.}\\[4mm]

\noindent{\bbb{References}}

\end{document}